\documentclass[12pt]{article}

\usepackage{amsmath,amsthm,amssymb}

\theoremstyle{plain}
\newtheorem{thm}{Theorem}
\newtheorem{prop}{Proposition}

\newtheorem{cor}{Corollary}
\newtheorem{lem}{Lemma}[section]

\def\R{\mathbb{R}}
\def\Z{\mathbb{Z}}

\def\v2{\vskip2mm}
\def\n{\noindent}

\def\({(\!(}
\def\){)\!)}
\def\e{\varepsilon}
\def\de{\delta}

\def\k{\kappa}
\def\la{\lambda}

\def\th{\theta}

\def\De{\Delta}

\def\La{\Lambda}
\def\Om{\Omega}

\def\pf{{\it Proof.}}
\def\v2{\vskip2mm}
\def\n{\noindent}

\def\be{{\bf e}}
\def\0{{\bf 0}}

\def\tst12{{\textstyle \frac12}}

\def\n{\noindent}
\def\beq{\begin{eqnarray*}}
\def\eeq{\end{eqnarray*}}

\def\beqn{\begin{equation}}
\def\eeqn{\end{equation}}

\begin{document}

\begin{center}
{\bf Asymptotic behaviour of  a random walk  killed on a finite set} \\
\vskip4mm
{K\^ohei UCHIYAMA} \\
\vskip2mm
{Department of Mathematics, Tokyo Institute of Technology} \\
{Oh-okayama, Meguro Tokyo 152-8551\\
e-mail: \,uchiyama@math.titech.ac.jp}
\end{center}

\vskip8mm
\n
{\it running head}:  random walk killed on a finite set

\vskip2mm
\n
{\it key words}: two-dimensional random walk; 
  exterior domain; transition probability; overshoot estimates 
\vskip2mm
\n

\vskip6mm

\begin{abstract}
We study  asymptotic behavior, for large time $n$, of the   transition probability of a  two-dimensional random walk killed when entering into a  non-empty finite subset $A$.  We show that it behaves like  $4 \tilde u_A(x) \tilde u_{-A}(-y) (\lg n)^{-2} p^n(y- x)$ for large  $n$, uniformly  in the parabolic regime $|x|\vee |y| =O(\sqrt n)$, where  $p^n(y-x)$ is the transition kernel of the random walk  (without killing)  and $\tilde u_A$ is the unique harmonic function in the \lq exterior of $A$' satisfying the boundary condition  $\tilde u_A(x) \sim \lg |x|$ at infinity.
 \end{abstract}
\vskip6mm

\section { Introduction and main results}

Let $S_n=S_0+X_1+\cdots+X_n$, $n=1,2,\ldots$ be a  random walk on the $d$-dimensional square lattice $\Z^d$, $d\geq 2$, defined on some probability  space $(\Om, {\cal F}, P)$. Here the  increments $X_j$ are i.i.d. random variables  taking values in $\Z^d$; the initial state $S_0$ may be any random variable specified according to the occasion.  As  usual the conditional law given  $S_0=x$  of the walk  $(S_n)$  is denoted by  $P_x$ and the  expectation under it by $E_x$.  
Let $A$ be a  non-empty finite subset of $\Z^d$. We are concerned with the asymptotic behavior of
$$p_A^n(x,y) = P_x[ S_k\notin A\,\,\mbox{for}\,\, k=1,\ldots, n\,\, \mbox{and}\,\, S_n =y],$$
the transition probability matrix of the walk killed on hitting $A$.  In  the classical paper \cite{K} H. Kesten proved among others that if  $d=2$ and   the walk is   recurrent  and (temporally)    aperiodic,  then for each  $x\in \Z^2$ and $y \in \Z^2\setminus  A$, as $n\to\infty$ 
\beqn\label{r1}
p^n_A(x,y) = u_A(x)u_{-A}(-y)\sum_{\xi\in A} q_A(\xi, n)(1+o(1)).
\eeqn
Here $q_A(\xi, n)$ is the 
$P_\xi$-probability of the walk returning to $A$  at time $n$ for the first time and $u_A$ is a unique harmonic function for the killed walk (i.e., $E_x[u_A(S_1); S_1\not \in A] =u_A(x)$ for all $x\in \Z^2$) that satisfies   $\sum_{\xi\in A}u_A(\xi)=1$ (the formulation is slightly modified from \cite{K} in which the function $P_x[ S_k\notin A\,\,\mbox{for}\,\, k=1,\ldots, n-1\,\, \mbox{and}\,\, S_n =y]$ is considered  instead of  $p^n_A$). 
 In this paper we improve the above asymptotic formula under  existence of the  finite second moments, so that  it is  valid  uniformly within the parabolic region  $|x|\vee |y| <M \sqrt n$ for each $M\geq 1$. Our approach  is different from that of \cite{K} and does not depend on the result of \cite{K}. 
  The proof in \cite{K} is done by  compactness arguments, namely  by showing      the convergence of the ratio $p^n_A(x,y)/ \sum_{\xi\in A} q_A(\xi, n)$ as  $n\to\infty$ along subsequences and  identifying    the limit. 
  In our approach, suggested by   intuitive idea  of how  random walk paths must behave to 
  contribute  to the transition probability,  we directly compute  the limit by using  an  estimate of the probability of the walk making  a large excursion  without hitting  $A$.
 For the estimation of such a probability  the potential function of the random walk with  its fundamental properties established by Spitzer \cite{Sh} plays a substantial  role as in \cite{K}.

  For one dimensional recurrent and  aperiodic walk  Kesten \cite{K} proved that if  $E X^2 =\infty$, then (\ref{r1}) is valid, and
if $E X^2=\sigma^2<\infty$  and $A=\{0\}$, then for each  $x, y\neq 0$
$$p^n_A(x,y) = [a(x)a(-y) + \sigma^{-4}xy]q_{\{0\}}(0, n)(1+o(1)),$$
where $a(x)$ is a potential function of the walk.
The problem of extending the latter result to a general  $A$ consisting of  more than one points is studied in a separate paper \cite{U1dm}.  In the higher   dimensional case $d\geq 3$,  where Kesten \cite{K} also gives a definite result similar to  (\ref{r1})  in a quite general framework, the uniform estimate is readily obtained if the existence of  the second moments is assumed as we shall briefly mention at the end of this section.

The corresponding problem for  two-dimensional Brownian motion is dealt with by the present author \cite{Ubdsty}. Although the strategy of the proof is the same for the both processes,    the Brownian case is basically  simpler and for it we can obtain quite 
detailed estimates valid uniformly beyond  parabolic region, whereas  adaptation of the  proof to  the random walk case, requiring  us to manage the overshoots to obtain relevant estimates, is nontrivial even though the space-time parameters  are restricted to a parabolic region. 

Throughout this paper we suppose that the random walk  $S_n$ is irreducible
 (i.e., for every $x\in \Z^d$,     $P_0[ S_n=x]>0$ for some $n>0$)
and that  
\beqn\label{mom}
EX=0  \quad  \mbox{and} \quad  E|X|^{2}<\infty. 
\eeqn
Here $X$ is a random variable having the same law as $X_1$ and $|\cdot |$ the usual Euclidian norm.
Let 
$p^n(x)= P_0[S_n=x]$
 so that $p^0(x)=\de(0,x)$ (Kronecker's delta kernel) and 
 $$P_x[S_n=y] =p^n(y-x).$$ 

For a non-empty set $B\subset \Z^d$,  $\sigma_B$ (resp. $\tau_B$)  denotes the first time when  $S_n$   enters into (resp. exits from) $B$:
$$\sigma_B = \inf\{n\geq 1: S_n \in B\},  \quad\quad  \tau_B = \inf\{n\geq 1: S_n \notin B\}.$$
(We shall sometimes write $\sigma(B)$ for $\sigma_B$ for typographical reason and similarly for $\tau(B)$.) 
By means of $\sigma_B$  the transition probability of $S_n$ killed on entering $B$ is written  as
$$p_B^n(x,y) = P_x[S_n =y, \sigma_B >n],\quad n=0,1, 2,\ldots;$$
in particular
$$p_B^0(x,y)= \de(x,y)\quad\mbox{and}\quad p^1_B(x,y) =p(y-x)(1-\chi_B(y))$$
 for all $x$, $y$, where $\chi_B(y) =1$ or $0$ according as $y\in B$ or $\notin B$.
Note that    $p_B^n(x,y) =0$ whenever $y\in B, n\geq 1$.

Let $A\subset \Z^d$  be a non-empty finite set such that the walk killed on $A$ is irreducible, so that 
\beqn\label{ass1}
\forall x \notin A, \forall  y\notin A,\,  \exists n\geq 1,  \quad p_A^n(x,y) >0,
\eeqn
which is supposed throughout the sequel. This imposes no essential restriction (see Remark 1(a)  after Theorem \ref{thm1}). 

Let $d=2$. Denote  the Green function associated with $p_A^n$ by $g_A(x,y)$: 
$$g_A(x,y) = \sum_{n=0}^\infty p_A^n(x,y),\quad  x, y\in \Z^2$$
 and put
\beqn\label{def_u}
u_A(x)= \lim_{|y|\to\infty} g_A(x, y),\quad x\in \Z^2.
\eeqn
(Cf. \cite[Theorem 14.3]{S} for the existence of the limit.)   $u_A(x)$ may be interpreted as the expected number  of visits  to $x$ made by  the dual  (or time-reversed)  random walk \lq starting at infinity'  up to (inclusively) the time of  first entrance into   $A$.
Let $Q$ be the covariance matrix of $X$, namely the  $2\times 2$ matrix whose quadratic form equals $E(X\cdot \th)^2$, $\th\in {\mathbb R}^2$,  and put
$$\k = \pi\sqrt{ \det Q}. $$
As we shall see in Section 2.3  it holds that as $|x|\to\infty$
\beqn\label{propty_u}
\k \;\! u_A(x) \,\sim\, \log |x|.
\eeqn
 Note that for $x\in A$,  $u_A(x)>0$ if and only if $p(y-x)>0$ for some $y\notin A$ (under condition (\ref{ass1})). 
\v2
\begin{thm}\label{thm1} Let $d=2$ and $A$ be a finite subset of $\Z^2$ that is non-empty and satisfies (\ref{ass1}). Then,  for each $M\geq 1$,  uniformly for $x\in \Z^2$ and  $y \in  \Z^2\setminus A$ subject to the constraint $|x|\vee |y| <M\sqrt n$, as  $n\to\infty$ 
\beqn\label{eq1.1}
p_A^n(x,y) = \frac{4\k^2 u_A(x)u_{-A}(-y)}{( \lg n)^2}p^n(y- x)(1+o(1)). 
\eeqn
\end{thm}
\v2
\n
{\sc Remark 1.} 
(a)\, If condition (\ref{ass1}) is not assumed, it  may possibly   occur  that $u_A(x)u_{-A}(-y)$ $ =0$ and $p^n_A(x,y)>0$ for some $x\in \Z^2$ and $y\notin A$ and for infinitely many $n$; thus  (\ref{eq1.1}) is not always true. However, if $u_A(x)u_{-A}(-y) =0$, then the random walk paths that connect $x$ and $y$ and avoid $A\setminus \{x,y\}$ are all confined in any disc that contains  $A$, hence  $p^n_A(x,y)$
approaches zero exponentially fast and we may consider only  the case  $u_A(x)u_{-A}(-y) >0$, or,  what amounts to substantially  
the same thing, augment $A$ by adding all $x$ with $u_A(x) u_{-A}(-x)=0$ so that (\ref{ass1}) is satisfied by the resulting set.

(b)\,  Under the stronger moment condition  $E[|X|^{2+\de}]<\infty$ (for some $\de>0$) the error term $o(1)$ can be replaced by $o( [\lg \lg n]/\lg n)$ in (\ref{eq1.1}). Although we shall not give full proof of it, some estimates  needed for it will be provided.

(c)\, Our definition of $g_A$  is not standard.  In \cite{S} and \cite{LL} the Green 
function for the walk killed on $A$ is given by  $G_A(x,y)= \sum_{n=0}^\infty P_x[S_k\notin A (k=0,\ldots,n), S_n=y] $ so that $G_A(x,y)=0$ if either $x\in A$ or $y\in A$, but still    $g_A(x,y) =G_A(x,y)$ whenever $x\notin A$. If  the hitting distribution  $H_A(x,y)$ is defined to be equal to  $ P_x[S_{\sigma(A)} =y]$ if $x\notin A$ and $\de(x,y) $ if $x\in A$, then  its dual $\widehat H_A(x,y)$ say, i.e., the corresponding distribution for the dual process $\widehat S_n:=\widehat S_0-X_1 -\cdots -X_n$, being equal to  $H_{-A}(-x,-y)$, we have
$$g_A(x,y) =G_A(x,y) +H_{-A}(-y,-x).$$ 
By means of this $g_A$ we shall have a neat expression of $u_A(x)$ (see Lemma \ref{lem2.1}).

\v2
Taking limit  in  $g_A(x,y) = \de(x,y)+  \sum_{z\notin A} p(z-x)g_A(z,y)$
we deduce   that
$u_A$ is harmonic  for the walk killed on $A$ in the sense that
\beqn\label{id} u_A(x) = E_x[u_A(S_1); S_1\notin A]\quad \mbox{for all}\quad x\in \Z^2.
\eeqn
For  $\xi\in A$, in particular, we have $u_{-A}(-\xi) = \sum_{y\notin A} u_{-A}(-y)p(\xi-y).$ 
Keeping this in mind  we substitute the expression given in  Theorem \ref{thm1} for $p_A^{n-1}(x,y)$   in  the identity
$$P_x[ \sigma_A =n, S_{n} =\xi] =\sum_{y\notin A} p_A^{n-1}(x,y)p(\xi-y),$$
and notice  that in view of a local limit theorem  
 $p^{n-1}(y-x)$ is well approximated by $p^n(\xi-x)$ for any sufficiently  large $n$ and any $y$ with $p^{n-1}(y-x)p(\xi-y)>0$ uniformly under the constraints $|y| =o(\sqrt n)$  and $|x|<M\sqrt n$, which leads to  the following 

\begin{cor}\label{cor1} Let $d=2$ and $A$ be a finite subset of $\Z^2$ that is non-empty and satisfies (\ref{ass1}). Then,  for each $M>1$ and  for $\xi\in A$,  uniformly for $x\in \Z^2$ with 
$|x|< M\sqrt n$,  as $n\to\infty$ 
\beqn\label{2}
P_x[ \sigma_A =n, S_{n} =\xi] = \frac{4\k^2u_A(x) u_{-A}(-\xi)}{(\lg n)^2}p^n(\xi-x) (1+o(1)).
 \eeqn 
\end{cor}
\v2

  The function   {\it $u_A$ restricted on  $A$  may be regarded as the hitting distribution of  $A$ for the dual  walk started at infinity} (since the walk visits  $A$ exactly once in the interval $[0,\sigma_A]$; cf. also Remark 1 (c)). It is noted that this fact may be expressed as 
$$u_{-A}(-\xi)= \lim_{|y|\to\infty}  P_y[ S_{\sigma_{A}}=\xi], \quad \xi\in A.$$

It  follows that $\sum_{\xi\in A} u_A(\xi) =1$ and  we deduce from (\ref{2}) that   {\it if the walk is  aperiodic, then}
$$P_x[ \sigma_A =n] = \frac{4\k^2u_A(x)}{(\lg n)^2}p^n(-x) (1+o(1))$$
and
$$\sum_{\xi\in A}P_\xi[ \sigma_A =n] = \sum_{\xi\in A}q_A(\xi,n) =\frac{4\k^2}{(\lg n)^2}p^n(0) (1+o(1));$$
in particular the formula (\ref{eq1.1}) conforms to (\ref{r1}) for each $x, y$ fixed.

The following proposition, crucial in our proof of Theorem \ref{thm1}, may explain why the function $u_A(x)$ comes into formula (\ref{eq1.1}).  Let  $U(R)$, $R>0$  denote the disc of radius  $R$:
$$U(R)= \{x\in \Z^2: |x|<R\}.$$
\begin{prop}\label{prop2.1}\, Uniformly for $x \in U(R)$, as $R\to\infty$
$$P_x[\tau_{U(R)} < \sigma_A] =\frac{ \k \;\! u_A(x)}{\lg R} (1+ o(1));$$ 
and if $E[X^2\lg |X|] <\infty$, then the error term  $o(1)$ can be replaced by $o(1/\lg R)$.  
\end{prop}

For the description of the strategy of the proof of Theorem \ref{thm1} the  readers are referred to   \cite{Ubdsty}: in the latter half of its first section   the skeleton of  proof to the corresponding
 result for Brownian motion  is given,   by which  the role of Proposition \ref{prop2.1} may be fully  understood.

Proof of Proposition \ref{prop2.1} is given at the end of Section 2, where some known results and immediate 
 consequences  of them are stated, overshoots estimates are discussed, and the function $u_A$  is defined in a way  apparently    different      from (\ref{def_u}).   Theorem \ref{thm1} is proved in Section 3 after two propositions    are 
 shown. 

\v2
We conclude this section by  making a short mention of the higher dimensional case. Let  $d\geq 3$ and define
$$u_A(x) = P_x[ \sigma_A =\infty].$$
As an  analogue of Theorem \ref{thm1} we can then verify   that uniformly for $x\in \Z^d$ and $y\in \Z^d \setminus A$ satisfying $|x|\vee |y|< M\sqrt n$, as $n\to\infty$
\beqn\label{d_3}
p_A^n(x,y) = u_A(x)u_{-A}(-y)p^n(y-x)(1+o(1)).  
\eeqn
 In view of the bound
 $$  \sum_{k=0}^\infty p^k(x) =o(1) \quad  (|x|\to\infty)$$
 as well as the  trivial facts: $P_x[\sigma_{\{0\}} =n] \leq p^n(-x)$;  $u_A(x)\to 1$ ($|x|\to\infty$); and
 $$P_x[ \sigma_{U(R)} <\sigma_A] \,\to\, u_A(x)\quad (R\to\infty)$$
the proof of Theorem \ref{thm1} is easily adapted for  verification of the formula (\ref{d_3}) with much 
simplification.

\section{Overshoot estimates and proof of Proposition \ref{prop2.1}} 
\subsection{Preliminary results}

Here we collect known results or its refined version that are used in the succeeding sections.
Let $d=2$. Put
$$a^\dagger(x) = \de(x,0) + a(x),$$ 
where $\de(x,y) = 0$ or $1$
according as  $x\neq y$ or $x=y$.
We recall that for all  $x, y$
$$ a(x) := \sum_{n=0}^\infty [p^n(0) -p^n(-x)] \geq 0, $$
\beqn\label{p20}
 E_x[a(S_1-y)] =  a^\dagger(x-y),
\eeqn
and
\beqn\label{p2}
g_{\{0\}}(x,y) =  \de(x,0)  + a(x)+ a(-y) - a(x-y);
\eeqn
and that
as $|x|\to \infty$
\beqn\label{p1}
a(x) = \frac1{\pi \sqrt{\det Q}} \lg |x| + o(\lg x), 
\eeqn
\beqn\label{p10}
a(x+y) -a(x) \to 0 \quad \mbox{for each}  \,\,\,  y
\eeqn
(cf. \cite[Theorem 4.4.6]{LL} for (\ref{p1}) and \cite[Propositions  11.6 and 12.2]{S}  for (\ref{p2}) and  (\ref{p10})). By (\ref{p2}) one can readily verify that if $\xi\in A$, the  process $M_n:=a(S_{n\wedge \sigma_A} -\xi)$ is a martingale relative to the filtration of  the stopped process $S_{n\wedge \sigma_A}$ under $P_x$ for $x\neq \xi$.

 We write  $x^2$ for  $|x|^2$.  In  Lemmas \ref{lem2.2} and \ref{lem2.4}  below  $\nu$ denotes   the period of the walk $S$, i.e.,  $\nu$ is  the {\rm g.c.d.} of  $\{n:p^n(0)>0\}$.  
\v2
\begin{lem}\label{lem2.2}\,Let $d=2$.   Then, {\rm (i)} uniformly for  $x \in {\mathbb Z}^2$ with $p^n(-x)>0$, as $n\to\infty$
\beqn\label{2.30}
P_{x}[\sigma_{\{0\}} =n] = \frac{2\nu \k \;\! a^\dagger(x)}{n(\lg n)^2}(1+o(1))  +  O\bigg(\frac{x^2\vee 1}{n^2 \lg n}\bigg),
\eeqn
and {\rm (ii)} as $n\wedge |x| \to\infty$
\beqn\label{2.3}
P_{x}[\sigma_{\{0\}} =n] = \frac{4\k \lg|x|}{(\lg n)^2}p^n(-x)  +  o\bigg(\frac1{n \lg n}\Big( 1\wedge \frac{\sqrt n}{|x|} \, \Big)\bigg).
\eeqn
\end{lem}

\v2\n
\pf\, These are proved in \cite{Ufh}: the first formula is  a reduced version of  Theorem 1.3, the second  is  the first half  of Theorem 1.4 (see also  the next lemma). \qed
\v2

\v2
The next result taken from \cite[Corollary 6]{Um_add}  is a version of  local central limit theorems for the random walk on $\Z^d$, $d\geq 1$. (Corollary 6 of \cite{Um_add}, stated for  aperiodic  case, is adapted to
 general case in an obvious way.)  Put
$$p^{(d)}_t(x) =\frac{1}{(2\pi t)^{d/2}} e^{- x^2/2t}$$
and
$$\sigma= (\det Q)^{1/2d}, \quad \tilde x= \sigma \sqrt{Q^{-1}} x\quad \mbox{and}\quad  \|x\| = |\tilde x|= \sigma\sqrt{x\cdot Q^{-1} x}.$$
(The factor $\sigma$ is put to make  $\|x\|$ agree with $|x|$ when  $Q$ is isotropic.)
\begin{lem}\label{lem2.4}~ Suppose  $E[\, |X|^{2+\de}] <\infty$  for a  constant $0\leq \de<1$.  Then  for $n\geq1$ and $x$ with $p^n(x)>0$,
\begin{eqnarray}\label{eq:LLT40}
p^{n}(x)
= \nu p_{\sigma^2 n}^{(d)}(\tilde x)
+\frac1{n^{d/2}}\times o\Big(\frac{n }{(\sqrt n\vee |x|)^{2+\de}}\Big),\nonumber
\end{eqnarray}
as $n+|x|\to \infty$ (for each $d=1, 2,\ldots$).
\end{lem}

\begin{lem}\label{lem2.3}\, Let $d=2$ and $n\geq 2$.  Then as  $|x|\wedge n\to\infty$
\beqn\label{n/lg}
P_{x}[\sigma_{\{0\}} <n]   
=  \frac{1}{\lg n}\int_{\tilde x^2/n}^\infty e^{- u/2\sigma^2}\frac{du}{u}+  
 o\bigg( 1\wedge \frac{ \sqrt n}{|x|\lg n} \, \bigg); 
 \eeqn
 in particular $P_{x}[\sigma_{\{0\}} <n]  \to 0$ if  $\, \liminf [(\lg x^2)/\lg n] \geq 1 $.
\end{lem}
\v2\n
\pf\, For simplicity  suppose that  the walk is  aperiodic. First suppose  $|x|>\sqrt n/\lg n$. Using   (\ref{2.3})  for $\lg n< k < n$  and   $P_{x}[\sigma_{\{0\}} =k]  \leq p^k(-x) $ for   $k\leq \lg n$ we compute the sum of  $P_x[\sigma_{\{0\}}=k]$ over  $k<n$. By Lemma \ref{lem2.4} we have  $2\k p^k(-x)= k^{-1}e^{- \la \tilde x^2/k} +o(1/(k\vee x^2))$, where $\la = 1/2\sigma^2$, and with the sum $\sum_{k= \lg n}^n   (\lg k)^{-2}e^{- \la \tilde x^2/k}/k$ computed in a usual manner  we observe 
\begin{eqnarray}\label{dist0}
 P_{x}[\sigma_{\{0\}} <n] 
& =&  o\bigg(\frac{\lg n}{x^2}\bigg)\, +\, 2 \lg |x|
 \int_{\tilde x^2/ n}^{\tilde x^2/\lg n} \frac{e^{- \la u}u^{-1}du}{[\lg (\tilde x^2/u)]^2}(1+o(1))  \nonumber\\
&&
+ \sum_{\lg n<k < n} \bigg(  \frac{\lg |x|}{(\lg k)^2(k\vee x^2)} + \frac{1}{|x|\sqrt k \lg  k}  \bigg) \times o(1)\nonumber\\
&=& \frac{2\lg |x|}{(\lg n)^2}\int_{\tilde x^2/n}^\infty e^{-\la u}\frac{du}{u} + o\Big(\frac{\sqrt n}{|x| \lg n}\Big).
\end{eqnarray}

Next suppose $|x|\leq \sqrt n/\lg n$ and   effect a similar    summation over $k\geq n$  of the right side of  (\ref{2.30})  (with  $k$ replacing  $n$)  and  subtract the resulting sum  from unity, and  then, recalling (\ref{p1}) you   find that 
\begin{eqnarray}\label{dist}
P_x[\sigma_{\{0\}} < n] &=& 1- \frac{2\k \;\! a^\dagger(x)}{\lg n}(1+o(1))  +  O\bigg(\frac{x^2}{n \lg n}\bigg)  \\
&=& \frac{0\vee \lg(n/ x^2)}{\lg n}+  o(1). \nonumber
\end{eqnarray}
By a little reflection these two relations  show  (\ref{n/lg}) as desired.  \qed

\v2\n
{\sc Remark 2.} \, Formula (\ref{n/lg}) does not identify the leading order of $P_x[\sigma_{\{0\}} < n] $ even within $|x| =O(\sqrt n)$: eg., if $|x|\sim \sqrt n/\lg n$, the error term reduces to $o(1)$ while the leading term is  $O((\lg\lg n)/\lg n)$.   However,  if $E[X^2\lg |X|]<\infty$, then $ a(x) = \k^{-1}\lg \|x\| + c^* +o(1)$ as $|x|\to\infty$  for some constant $c^*$ \cite[Theorem 1]{Ugf} and   from (\ref{dist})
  we have
$$P_x[\sigma_{\{0\}} < n] = \frac{ \lg(n/\tilde x^2)}{\lg n}+ \frac{2\k c^*}{\lg n}(1+ o(1)) + O\bigg(\frac{x^2}{n \lg n}\bigg)$$
as $|x|\wedge n \to\infty$, which combined with (\ref{dist0}) (valid for $|x|\geq \sqrt n /\lg n$) provides  an explicit  asymptotic form of
the probability within the regime  $|x| =O(\sqrt n)$.

\v2
Let $(\eta_n)_{n=1}^\infty$ be a sequence of  real i.i.d. random variables of mean zero  with a finite and positive variance,  and consider the random walk $Y_n = \eta_1+\cdots+ \eta_n$ on $\R$.    For convenience of   later citation we record the following two lemmas. Although it may  fall into a body of common knowledge,  proofs are given for completeness.  
\begin{lem}\label{lem2.5} Let $Y_n$  be as above and  $\sigma^Y_{B}$ denote the first entrance time  of the walk $(Y_n)$ into $B$. Then, as $R\to\infty$ 
$$\frac1{R} E[ Y_{\sigma^Y_{[R,\infty)}} -R\, ] \, \longrightarrow\, 0.$$
\end{lem}
\v2\n
\pf\,  Let $Z$ be  the first strict ascending ladder variable in  $(Y_n$), i.e.,   $Z = Y_{\sigma^Y_{(0,-\infty)}}$ and  $V(dx)$  the renewal  measure of the ascending ladder process:  $V(dx) =\sum_{n=1}^\infty P[Z_n \in dx]$, where $Z_n$ is the  $n$-th record value  of   $Y$. Then
$$E[ Y_{\sigma^Y_{[R,\infty)}} -R\,] = \int_{(0, R)} V(dy) \int_{[R-y, \infty)} z P[Z\in dz]. $$
On knowing $E Z<\infty$ \cite[Theorem 18.5.1]{F} and  $V((0,R]) \leq CR$, this is  $o(R)$ as readily ascertained. \qed

\begin{lem}\label{lem2.40} \, Let $Y_n$  be as above and $\sigma^2 :=E Y_1^2$.  Then there exists a universal constant $0<c_0<1$ such that for $N > 0$ and for all sufficiently large  $R$,
   $$ \sup_{x\in \R: |x|<R} P[\,|Y_k+x| < R \,\; \mbox{ for}\,\; k \leq N\,] \leq  c_0^{-1+ \sigma^2N/R^2}.$$
\end{lem}
\v2\n
\pf\,   From  central limit theorem it follows that 
   $$\limsup_{R\to\infty} \sup_{|x|< R} P[\,|Y_k+x| < R \,\; \mbox{ for}\,\; k\leq R^2/\sigma^2] 
   \leq  \lim_{R\to\infty} P[\, |Y_{\lfloor R^2/\sigma^2\rfloor}| \leq  R\,] =: c,$$
 where $c = \int_{-1}^1 \frac{e^{- u^2/2} du}{\sqrt {2\pi}}$ .  Thus,   taking any $c_0$ from $(c,1)$ we find
  that  if $n = \lfloor \sigma^2 N/R^2\rfloor$, 
$$ P[\,|Y_k+x| < R \,\;\mbox{for} \,\;k\leq  N] \leq c_0^n \qquad (|x|<R)$$
for all sufficiently large $R$. \qed

\subsection{Overshoot estimates}
\begin{prop}\label{prop4.1}\, Uniformly for $x\in U(R)$, as $R\to\infty$
$$ \frac1{R}E_x [ \,|S_{\tau_{U(R)}}|-R \,] \,\longrightarrow \,  0.$$
\end{prop}
\v2\n
\pf\,  It suffices to show that  for each $\e>0$,
\beqn\label{3.1}
\frac1{R}E_x\Big[\, |S_{\tau_{U(R)}}| -R \, ;  |S_{\tau_{U(R)}}|  > (1+\e )R \Big] \,\longrightarrow \,  0.
\eeqn
For the proof we divide the plane  by several rays emanating from the origin to reduce the problem to a one-dimensional one. Taking a positive integer $N$ so that
\beqn\label{1+2}
(1+\e)\cos (\pi/N) \geq 1, 
\eeqn
we  put
$$W_1= \{y=(y_1,y_2) \in \Z^2 : y_1 \geq  (|y_2|\tan {\textstyle \frac\pi{N}}) \vee R \, \},$$
the intersection of the half plane  $y_1\geq R$ and the infinite sector  held by the  two rays 
$y_1=\pm y_2 \tan \frac \pi{N}, y_1 >0$,
and let $W_k$  $(k=2,\ldots,  N)$ be the rotation of $W_1$ by $2k\pi/N$ 
counterclockwise. Then, by virtue of  (\ref{1+2})  the event
$\{|S_{\tau_{U(R)}}| >  (1+\e) R\}$ occurs only if  one of  the $N$ events
$$\La_k:= \{ S_{\tau_{U(R)}}  \in W_k \}, \quad k=1,2, \ldots, N$$
occurs, hence the expectation to be estimated is dominated by 
$\sum E_x[ \, |S_{\tau_{U(R)}}| - R \, ; \La_k]$,  
of which we  consider the term with $k=1$ for convenience of description. Let  $Y_n$ be
 the first component of $S_n$: $Y_n = \be \cdot S_n$.  Noting that the occurrence of $\La_1$  entails 
 $Y_{\tau_{U(R)}} \geq R$  we  see   
\[
 E_x[\, |S_{\tau_{U(R)}}| -R\, ; \La_1]  \leq  2 E_x[  Y_{\sigma^Y_{[R,\infty)}} -R\, ;  \tau_{U(R)} = \sigma^Y_{[R,\infty)}\, ]  =o(R).
 \]
 Here the last equality follows from  Lemma \ref{lem2.5}, for the expectation on the middle member is at most
  $E_0[ Y_{\sigma^Y_{[R',\infty)}} -R'\,]$ with  $R'= R- \be\cdot x \leq 2R$. 
By the same argument we obtain  $ E_x[\, |S_{\tau_{U(R)}}| -R\, ; \La_k] =o(R)$ for each $k$, 
Thus   the  relation (\ref{3.1})  is verified. \qed

\v2
By Markov's inequality there follows immediately  from Proposition \ref{prop4.1} the following
\begin{cor}\label{cor4.1}\, Uniformly for $x\in U(R)$ and for $R'>R$, as $R\to\infty$
\beqn \label{(ii)}
P_x[\, |S_{\tau_{U(R)}}| >  R'\, ] = o\Big( \frac{R}{R'-R}\Big).
\eeqn 
\end{cor}

\begin{lem}\label{lem2.50}\, Uniformly for $x\in U(R)$, as $R\to\infty$
\beqn\label{eq_lem2.5}
E_x\Big[ \lg(|S_{\tau_{U(R)}}| / R) \,; \tau_{U(R)} < \sigma_A\Big] = P_x[ \tau_{U(\sqrt R)} < \sigma_A]\times o(1)
+ o(R^{-1/2}).
\eeqn
\end{lem}
\v2\n
\pf\,  Write $\De = \lg (|S_{\tau_{U(R)}}| /  R) $.   
 By the inequality  $\lg (u/R)\leq  (u-R)/R$\,   ($u>R$) 
 \[ 
 E_x[\, \De \,; \tau_{U(R)} < \sigma_A] 
   \leq  \frac{1}{R} E_x\Big[ |S_{\tau_{U(R)}}| -R\, ; \tau_{U(R)} < \sigma_A   \Big].
\]
The right side tending to zero according to  Proposition \ref{prop4.1},   (\ref{eq_lem2.5})   holds   if $\sqrt R\leq |x| <R$,  since then 
the probability on the right side of it  is bounded away from zero.

As for  $x \in U(\sqrt R)$  we decompose
$$E_x[\De ; \tau_{U(R)} < \sigma_A]  = E_x[ \De; \tau_{U(\sqrt R)}= \tau_{U(R)} < \sigma_A] +  E_x[ \De; \tau_{U(\sqrt R)}< \tau_{U(R)} < \sigma_A]. $$
The second term on the right side  is written as 
\[
 E_x[ E_{S_{\tau(U(\sqrt R))}}[\De; \tau_{U(R)} < \sigma_A];  \tau_{U(\sqrt R)}< \tau_{U(R)} \wedge  \sigma_A], 
 \]
which, owing to what we have observed above,  is evaluated to be
$ P_x[  \tau_{U(\sqrt R)}<  \sigma_A] \times o(1).$
On the other hand the first term  is dominated by
\beq
E_x[\De;  S_{\tau_{U(\sqrt R )}} \geq R] &\leq& 
2P_x[R\leq  S_{\tau_{U(\sqrt R )}} < R^2] + E_x[ \, \lg |S_{\tau_{U(\sqrt R )}}| ;  |S_{\tau_{U(\sqrt R )}}| \geq  R^2] \\
&\leq& o(R^{-1/2}) + R^{-1} E_x[ \,|S_{\tau_{U(\sqrt R )}}|^{1/2}\lg |S_{\tau_{U(\sqrt R )}}|\, ]\\
 &=& o(R^{-1/2}),
\eeq
where we have applied Corollary \ref{cor4.1} and Markov's inequality  for the second inequality and Proposition \ref{prop4.1} for the last relation. 
The proof of Lemma \ref{lem2.50} is complete. \qed

\begin{lem}\label{lem4.3}\, Uniformly for $x\in U(R)$, as $R\to\infty$
$$ E_x[\,|\k \;\! a(S_{\tau_{U(R)}}) - \lg R|\,; \tau_{U(R)} < \sigma_A]  = (1\vee \lg |x|)  \times  o(1).$$
If $E[ X^2\lg |X|] <\infty$, then the right side can be replaced by $ P_x[\tau_{U(\sqrt R)} < \sigma_A]\times O(1) + o(R^{-1/2})$.
\end{lem}
\v2\n
\pf\,  We may and do suppose $0\in A$, so that   $a$ is harmonic on $\Z^2\setminus A$. By  the optional stopping theorem applied to the positive martingale  $a(S_{n\wedge \sigma(A)})$ it follows that 
\[
 0\leq E_x[ a(S_{\tau_{U(R)}})  \,; \tau_{U(R)} < \sigma_A]  \leq a^\dagger(x)
\]
for all  $x$. (Note the case $x\in A$ is reduced to the case $x\notin A$.) 
In the decomposition
$$\k \;\! a(S_{\tau_{U(R)}}) - \lg R =\k \;\! a(S_{\tau_{U(R)}}) - \lg |S_{\tau_{U(R)}}| + \lg ( |S_{\tau_{U(R)}}|/ R), $$
the contribution of the first two terms on the right side is therefore  disposed of as follows:
\beq
&& E_x\Big[\, \Big|\k \;\! a(S_{\tau_{U(R)}}) - \lg |S_{\tau_{U(R)}}|\, \Big|    \,; \tau_{U(R)} < \sigma_A\Big] \\[2mm]
&&= \left\{ \begin{array} {ll}    a^\dagger(x)\times o(1) \quad  &\mbox{in general}, \\
 P_x[\tau_{U(R)} < \sigma_A]\times O(1) \quad &\mbox{if}\quad E[X^2 \lg |X|]<\infty,
 \end{array} \right.
\eeq
where we have also applied  the fact that $\k \;\! a(z) -\lg |z|$ is $o(a(z))$ ($|z|\to\infty$) in general and bounded  if $E[X^2 \lg |X|]<\infty$ (\cite[Theorem1]{Ugf}), whereas 
the required bound concerning the third term is provided by  Lemma \ref{lem2.50}. 
\qed 

\v2
In the next  lemma we  include the case when  the existence of higher moments are assumed  though not applied in this article.
\begin{lem}\label{lem4.4}\, Suppose $E[\,|X|^{2+\de}]<\infty$ for a  constant $\de\geq 0$. Then, uniformly for $R>1$ and  $x\in U(R)$, 
$$P_x[S_{\tau_{U(R)}} =y\, ,  \tau_{U(R)} < \sigma_A] = a^\dagger (x)\times o((|y|-R)^{-(2+\de)}) \quad ( |y|-R \to\infty),$$
 and
 \beq
  P_x\Big[ |S_{\tau_{U(R)}}|> R+ h \, ,  \tau_{U(R)} < \sigma_A\Big]  = a^\dagger (x) R^2\times 
o(h^{-(2+\de)})  \quad (h\to \infty).
 \eeq\end{lem}
\v2\n
\pf\, Take $R$ such that $A\subset U(R)$ and put  $D(R)= A\cup (\Z^2 \setminus U(R))$. Then  for $y\notin U(R)$,  
\beq
 P_x[S_{\tau_{U(R)}} = y \, ,  \tau_{U(R)} < \sigma_A] 
=\sum_{z\in U(R)\setminus A} g_{D(R)}( x,z)p(y-z).
\eeq
 Suppose  $0\in A$ for simplicity. Then,  $g_{D(R)}(x,z)\leq g_{\{0\}}(x,z)\leq g_{\{0\}}(x,x)\leq Ca^\dagger(x)$. Hence
$$ P_x[S_{\tau_{U(R)}} = y \, ,  \tau_{U(R)} < \sigma_A\,]   \leq Ca^\dagger (x) P[\,|y-X| <R] \leq Ca^\dagger (x) P[\, |X|>|y|- R\,] $$
($y\notin U(R)$) and an application of Markov's inequality yields  the first relation of the lemma. To verify the second relation we use the first inequality above and observe
$$\sum_{|y|> R+h} P[\,|y-X| <R]= \sum_{|z|<R}\sum_{|y|> R+h}p(y-z) = 4 R^2P[|X| > h],$$
which implies the required bound in view of  Markov's inequality. \qed

\subsection{The function $u_A$ and proof of Proposition \ref{prop2.1}}
\v2\v2
 We define a function $u_A(x)$ by
\beqn\label{p300}
u_A(x)= g_A(x,y) + a(x-y) - E_{x}[ a(S_{\sigma(A)}-y)], \quad  y\in \Z^2.
\eeqn
\begin{lem}\label{lem2.1}\,The right side of (\ref{p300}) is independent of $y\in \Z^2$ for all  $x\in \Z^2$. In particular for any $\xi_0\in A$
\beqn\label{p3}
u_A(x)= \de(x,\xi_0) + a(x-\xi_0) - E_{x}[ a(S_{\sigma(A)}-\xi_0)].
\eeqn
\end{lem}
\v2\n
\pf\,
With $x\in \Z^2$ fixed  put
$f(y):= g_A(x,y) + a(x-y) - E_x[ a(S_{\sigma(A)}-y)]$. We are to prove that  $f$ 
is dual-harmonic on $\Z^2$, namely $f(y)= \sum_z p(y-z) f(z)$ for all  $y\in \Z^2$, hence must be constant, for $f$ is bounded (cf. \cite[Proposition 6.3]{KSK}). Remember   that for all $x, y \in \Z^2$,
$$g_A(x,y) = \sum_{n=0}^\infty p_A^n(x,y), \quad p^0_A(x,y) =\de_{x,y};$$
in particular  $g_A(x,y) = \de_(x,y)$ for $y\in A$. 
We then observe  that if 
$\widehat P$ denotes the dual transition operator, i.e., $\widehat Pf(y)= \sum_z p(y-z)f(z)$, then for all $x\in \Z^2$,
\[   
\widehat P  a(x-\cdot)(y) = a(x-y)+ \de(x, y) 
\]
and
\[
\widehat P\{E_x[ a(S_{\sigma(A)}-\cdot)]\}(y) 
= \left\{ \begin{array} {ll}  E_x[ a(S_{\sigma(A)}-y)]  \quad &(y\notin A), \\[1mm]
 E_x[ a(S_{\sigma(A)}-y)] + P_x[ S_{\sigma_A}=y]  \quad & (y\in A),
 \end{array} \right.
\]
while
\[
\widehat P g_A(x,\cdot)(y) = \left\{ \begin{array} {ll}  g_A(x,y) - \de(x,y) \quad & (y\notin A),\\[1mm]
P_x[ S_{\sigma_A}=y] \quad & (y\in A).
\end{array} \right.
\]
From these there readily follows the identity $\widehat Pf =f$ as required. \qed
\v2\n
{\sc Remark 3.} \, The assertion of Lemma \ref{lem2.1} as well as the proof  given above is valid for all recurrent walks. For $x$ restricted on $A$,  (\ref{p300})  reduces to the dual of the formula of Proposition 30.1 in  \cite{S}, where the dual of $u_A(x)$, which  therein  is denoted by $\mu_A(x)$,  is defined as the limit of $\sum_{z\in \Z^2} p^n(z-y)P_z[S_{\sigma(A)} =x]$ as $n\to \infty$. On the other hand for $x\notin A$  an equivalent to (\ref{p300}) is verified in \cite[Proposition 4.6.3]{LL} in a different way when  $X$ is of  finite range and symmetric. 
Formula  (\ref{p300})    is suggested by a similar one  used by  Hunt \cite{Ht} to define  the (classical)  Green  function  in a plane region with pole at infinity.

\v2

 By virtue of   (\ref{p10}) passing to the limit in (\ref{p300}) yields
$$  \quad u_A(x)=\lim_{|y|\to\infty} g_A(x,y). $$
Thus  the present definition of  $u_A$  agrees with the one given in (\ref{def_u}).
By   (\ref{p1}) and (\ref{p3}) we also obtain that  $\k u_A(x) \sim \lg|x|$ ($|x|\to\infty$) as stated in (\ref{propty_u}).

Now we prove Proposition \ref{prop2.1}, which we state again
\v2

\v2\n
{\bf Proposition 1.}\, {\it Uniformly for $x \in U(R)$, as $R\to\infty$
$$P_x[\tau_{U(R)} < \sigma_A] = \frac{ \k \;\! u_A(x)}{\lg R} (1+ o(1));$$ 
and if $E[X^2\lg |X|] <\infty$, then the error term  $o(1)$ may be replaced by $O(1/\lg R)$.}  
\v2\n
\pf\, It suffices to prove the postulated formula for $x\notin A$,  the case $x\in A$  being reduced to it by virtue of  (\ref{id}). Suppose $x\notin A$ and 
 let $\xi\in A$ and  $D(R)=A\cup (\Z^2  \setminus U(R))$ (for $R$ large enough). We adapt   the proof in \cite{LL} (Proposition 6.4.7)  to the present situation of unbounded $X$ (cf. also  \cite{St} (Theorem 11.2.14)). By optional sampling theorem the process $M_n:=a(S_{n\wedge \sigma_{D(R)}}-\xi)$ is a martingale under $P_x$. It is easy to see that $M_n$ is $L^2$ bounded and we   apply the martingale
convergence theorem  to have 
$$a(x-\xi) =E_x[a(S_{\sigma_{D(R)}} -\xi)].$$
 Breaking this expectation 
according as the first visit to $D(R)$ occurs on $A$ or on $\Z^2\setminus U(R)$ yields
\beq
a(x-\xi) -E_x[a(S_{\sigma_A}-\xi)]  &=&E_x[ a(S_{\tau_{U(R)}} -\xi)\,;\, \tau_{U(R)} < \sigma_A]\\
&&  - E_x[ a(S_{\sigma_{A}} -\xi); \tau_{U(R)} < \sigma_A].
\eeq 
Recall that  by (\ref{p3})  the left side equals $u_A(x)$  ($x\notin A$).  Then,  putting 
$$c(x,R) = E_x[ a(S_{\sigma_{A}} -\xi) \,|\, \tau_{U(R)} < \sigma_A]$$
 and
\beqn\label{p5}
\th(x,R) =  E_x[ \;\! a(S_{\tau_{U(R)}} -\xi) -\k^{-1} \lg R\,;\, \tau_{U(R)} < \sigma_A],
\eeqn
  we rewrite  the identity above in the form
\beqn\label{p51}
u_A(x)  = \th(x,R) +  (\k^{-1}\lg R - c(x,R)) P_x[\tau_{U(R)} < \sigma_A].
\eeqn
Since  $c(x,R)$ is bounded and in view of Lemma \ref{lem4.3} and (\ref{propty_u})
 $\th(x,R) =  u_A(x) \times o(1)$ as  $R\to\infty$ uniformly for $x\in U(R)\setminus A$, on  dividing by $\k^{-1}\lg R$ we obtain   the formula  of the proposition. If $E[X^2\lg |X|]<\infty]$, by the second assertion  of Lemma \ref{lem4.3} and what is just proved we have $\th(x,R) = P_x[\tau_{U(R)} < \sigma_A]\times O(1) $, 
 and substitution of it  into (\ref{p51}) leads to the desired estimate of the error term. \qed

\section{Proof of Theorem \ref{thm1}}
Throughout this section let $d=2$.

\begin {prop}\label{lem3.1}   As   $|x|\wedge |y| \wedge n\to\infty$, 
 under the condition  $\liminf \, [ \lg (x^2\wedge y^2)]/\lg n\geq 1$,  
$$p_A^n(x,y) = p^n(y- x) + o\Big(\frac1{n \vee x^2\vee y^2}\Big).$$
\end{prop}
\v2\n
\pf\, Let $|x|\leq |y|$ and $p^n(y-x)>0$. 
We make decomposition
\beqn\label{100} 
p^n(y-x) - p^n_A(x,y) = \sum_{k=1}^n \sum_{\xi\in  A}  P_x[ \sigma_A =k, S_k =\xi]\, p^{n-k}(y-\xi),
\eeqn
which entails
\beqn\label{10} 
p^n(y-x) - p^n_A(x,y) \leq  \sum_{k=1}^n \sum_{\xi\in  A}  P_x[ \sigma_{\{\xi\}} =k]\, p^{n-k}(y-\xi).
\eeqn
Let $I_{[a,b]}$ denote the latter double sum   restricted to the interval $[a,b], 0\leq a<b\leq n$.

By  Lemma \ref{lem2.4}  (local  limit theorem) it follows that
for $ \xi \in A, k<n/2$,
\beqn\label{LLT}
p^{n-k}(y-\xi) \leq  p^n(y) +o(1/(n\vee y^{2})) =O(1/(n\vee y^2)), \quad 
\eeqn
 hence
\[
I_{[0,n/2]}  \leq \frac{C}{n\vee y^2} \sum_{\xi\in A}P_x[\sigma_{\{\xi\}}\leq n/2].
\]
 Applying Lemma \ref{lem2.3} to  
 $P_x[\sigma_{\{\xi\}}\leq n/2]$ 
 we deduce that
  \beqn\label{101}
I_{[0,n/2]}  \leq \frac{C'}{n\vee y^2}\bigg( \frac{1\vee \lg(n/x^2)}{\lg n}
 +o\Big(1\wedge \frac{\sqrt n}{|x|\lg n}\Big)\bigg).
\eeqn

For the other part  $I_{[n/2,n]}$ we  apply  (\ref{2.3}) and  observe  
\beqn\label{102}
I_{[n/2,n]} \leq  \frac{C}{n\lg n}\sum_{\xi\in A}\sum_{n/2 < k\leq n} p^{n-k}(y-\xi).
\eeqn
The inner sum $J_n := \sum_{n/2 < k\leq n} p^{n-k}(y-\xi)$ is estimated somewhat differently  depending on whether  $y^2>n/2$ or not. If $y^2\leq n/2$, then, on splitting the sum at $j:= n-k = \lfloor y^2\rfloor$ and  using Lemma \ref{lem2.4},
$$J_n =\sum_{j\leq y^2}  + \sum_{y^2<j \leq n/2}\leq \int_0^{y^2}\bigg(p^{(2)}_{\sigma^2 t}(\|y\|)+ o\Big(\frac1{y^2}\Big)\bigg)dt + \sum_{y^2\leq j\leq n/2}\frac1{j} =  O\Big(\lg({n}/{y^2})\Big),   $$
which we may write as
$$J_n \leq C'n p^n(y) \lg ({n}/{y^2}).$$
Similarly,  if $y^2> n/2$, then
\beqn\label{103}
J_n = \int_0^{n/2}p^{(2)}_{\sigma^2 t}(\|y\|) dt +  \sum_{1\leq j < n/2} o\bigg(\frac{1}{j\vee y^2}\bigg)  \leq  
 C\frac{n}{y^2}e^{-\|y\|^2/2\sigma^2n} + o\bigg(\frac{1}{y^2}\bigg) \times n. 
\eeqn
In view of (\ref{102}) these two bounds of $J_n$  show
\begin{eqnarray}\label{104}
I_{[n/2,n]}\leq\frac{C''\sharp A}{n\lg n}J_n =  p^n(y)\times O\bigg(\frac{1\vee\lg (n/y^2)}{\lg n}\bigg) + \frac1{\lg n}\times o\Big(\frac{1}{n\vee y^2}\Big).
\end{eqnarray}
By the \lq$\liminf$' condition we have $1\vee \lg(n/y^2) = o(\lg n)$, hence the right side above is $o(1/(n\vee y^2))$;
and similarly for the right side of  (\ref{101}).  This finishes the proof. \qed

\v2\n
{\sc Remark 4.} \,   The estimate of Proposition \ref{lem3.1}  may be improved under  a higher moment condition.   {\it  Suppose  $E|X|^{2+\de}<\infty$ for a constant $\de\geq 0$ and $|x|\leq |y|$.  Then, for any  $M>1$,  as  $|x|\wedge n\to\infty$ under the constraint $(1\vee|x|) |y| <Mn$}
 \beqn\label{Rem3}
 p_A^n(x,y) = p^n(y- x)\Bigg[1+ O\bigg(\frac{1\vee\lg (n/x^2)}{\lg n}\bigg)\Bigg] + o  \bigg(\frac{1}{( n\vee  y^2)^{1+\de/2}} \bigg).
 \eeqn

\v2\n
 The proof is carried out by examining the preceding one.  Owing to Lemma \ref{lem2.4},  $o(1/(n\vee y^{2}))$ may be replaced by  $o(( n \vee y^2)^{-1-\de/2})$ in 
(\ref{LLT}).
Taking Remark 2 (given to Lemma \ref{lem2.3})  into account, this leads to, instead of  (\ref{101}),
\beqn\label{105}
I_{[0,n/2]}  \leq \bigg[C p^n(y) + o\Big(\frac1{(n \vee  y^2)^{1+\de/2}}\Big) \bigg]
 \frac{1\vee \lg(n/x^2)}{\lg n}.
\eeqn
 As for $I_{[n/2,n]}$, using  (\ref{102}) as before  we infer that 
 on the right side of  (\ref{104}) $o((n\vee y^{2})^{-1})$  can be replaced by  $o(1/ ( n\vee y^2)^{1+ \de/2})$. Combined with (\ref{105})
 this shows   (\ref{Rem3}), since our supposition entails $|y-x|^2 - y^2= O(n)$, so that $p^n(y)$ can 
 be replaced by  $p^n(y-x)$. The details are omitted.

\v2
\begin {prop}\label{lem3.3} \,
As $n\wedge |y|\to\infty$  under the constraint 
  $\,\liminf {(\lg y^2)}/{\lg n} \geq 1,$
\beqn\label{eq3.3}
p_A^n(x,y) = \frac{2\k \;\! u_A(x)}{\lg n} \bigg[ p^n(y-x) + o\Big(\frac1{n \vee  y^2}\Big) \bigg]
\eeqn
 uniformly for $x \in U(\sqrt n)$.
\end{prop}
\v2\n
\pf\,  Suppose  $p^n(y-x)>0$. Define $R= R_{n}$  and $N= N_{n}$ by
$$R=   \frac{\sqrt n}{(\lg n)^2}   
\quad \mbox{and}\quad N= \lfloor (R\lg n)^2\rfloor,$$
so that
\beqn\label{R_N}
 (a)  \,\, \, \frac{N}{n} = O\Big(\frac{1}{(\lg n)^2}\Big) \quad\mbox{and}  \quad (b) \,\,\,   \frac{N+1}{R^2} \geq  (\lg n)^2.
\eeqn
The case $R\leq |x| < \sqrt n$ 
being  covered by Proposition \ref{lem3.1} because of (\ref{propty_u}),    in the sequel we suppose  that   $|x|< R$.

Define $\e(n,x,y)$ via
\begin{eqnarray} \label{13} 
p_A^n(x,y) &=& \sum_{k=1}^{N-1}\sum_{z\notin U(R)} P_x[ \tau_{U(R)} =k < \sigma_A, S_{\tau_{U(R)}}=z \,]\, p_A^{n-k}(z,y) \\
&& \, + \, \e(n,x,y).\nonumber 
\end{eqnarray} 
By Lemma \ref{lem2.40}  we have
\beqn\label{131}
P_x[\tau_{U(R)}\wedge \sigma_A \geq N] \leq c e^{-\la N/R^2}
\eeqn
with  a  constant $\la>0$ (depending  only on the walk)  and,  noting that for $z \in  U(R)$,  $p_A^{n-N}(z,y)\leq p^{n-N}(y-z)\leq c p^n(y-x) + o(n^{-1}\wedge y^{-2}) $, we observe 
\begin{eqnarray} \label{14}
\e(n,x,y) 
&=&\sum_{z\in U(R)\setminus A} P_x[ \tau_{U(R)}\wedge \sigma_A \geq N, S_N =z] \, p_A^{n-N}(z,y) \nonumber\\
&\leq & c'e^{-\la N/R^2}\Big[p^n(y-x) +o(n^{-1}\wedge y^{-2}) \Big],
\end{eqnarray}
hence $\e(n,x,y)$ is absorbed into the error term in (\ref{eq3.3})  because of (\ref{R_N}:(b)). 

To evaluate  the double sum in (\ref{13})  we take 
$$r=r_{y,n}=(\sqrt n \vee |y|)/\sqrt{\lg n},$$
 so that  
\beqn\label{R_N1}
 \frac{R^2}{r^2} = \frac{n}{(n\vee y^2)(\lg n)^3},
 \eeqn
and  first dispose of the part with the inner sum restricted to the outside of  $U(r)$.   Since for  $k\leq N$,
\[
p_A^{n-k}(z,y)  \leq C /n,
\]
we see that this part is
\begin{eqnarray} \label{ov_shoot2}
&& \sum_{k=1}^{N-1}\sum_{z\notin U(r)} P_x[ \tau_{U(R)} =k <\sigma_A, S_{\tau_{U(R)}}=z \,]\, p_A^{n-k}(z,y) \nonumber \\
&&\leq  \frac{C}{n} P_x[ \tau_{U(R)}  <\sigma_A, S_{\tau_{U(R)}} \not \in U(r) ],
\end{eqnarray}
hence negligible, for
by Lemma \ref{lem4.4} and (\ref{R_N1}) the probability on the right side is dominated by a constant multiple of $a^\dagger(x) R^2/r^2= O\Big(n/[(n\vee y^2)(\lg n)^{2}] \Big)$. 

For $k<N, \, z \in U(r) \setminus U(R)$ with $p^{n-k}(y-z) >0$, by Proposition \ref{lem3.1} we have
\[
p_A^{n-k}(z,y) = p^{n-k}(y-z) + o(1/(n\vee y^2));
\]
also we have  
\[
 p^{n-k}(y-z) = p^n(y-x) + o(1/(n\vee y^2)),
\]
where we have applied  $|z| = o(\sqrt n\vee |y|)$ and $|x|=O(\sqrt n)$. Thus the contribution to the double sum to be evaluated is written as
$$ P_x[\tau_{ U(R)} < N \wedge \sigma_A] \times [p^n(y-x) + o(1/(n\vee y^2)) ],$$
while,  owing to   Proposition \ref{prop2.1} and  (\ref{131}),
\begin{eqnarray}\label{24} 
P_x[\tau_{ U(R)} < N  \wedge \sigma_A] &=& P_x[\tau_{ U(R)} < \sigma_A] -  P_x[N\leq  \tau_{U(R)}  < \sigma_A]   \nonumber\\
&=& \frac{\k \;\! u_A(x)}{\lg R}(1+ o(1)).
\end{eqnarray}
Now, combining this with what we have observed by (\ref{ov_shoot2}),  we conclude that   the double sum in (\ref{13}) is written as 
\[
 \frac{\k \;\! u_A(x)}{\lg R}\, \bigg[p^n(y-x) + o\Big(\frac1{n\vee y^2} \Big) \bigg].
\]
  Since  $\lg R = \frac12 (\lg n)(1+o(1))$  as noted previously, this verifies the  lemma.  \qed
 
\v2\n
{\it Proof of Theorem \ref{thm1}. } \,  The relation to verify is
$$p_A^n(x,y) = \frac{4\k^2 u_A(x)u_{-A}(-y)}{(\lg n)^2}p^n(y- x)(1+o(1)).$$ 
In view of  Propositions \ref{lem3.1} and \ref{lem3.3} we may suppose  $|x|\vee |y| \leq \sqrt n /(\lg n)^2$.
 Consider   the double sum in (\ref{13}) (with the same $R$ and $N$ as therein),  of which we have observed that the contribution to it from  $|z|>\sqrt {n/\lg n}$ may be neglected, while 
by  Proposition \ref{lem3.3} we have  for $R\leq  |z|\leq \sqrt {n/\lg n}$ 
 $$p_A^{n-k}(z,y) = p^{n-k}_{-A}(-y,-z) = \frac{2\k \;\! u_{-A}(-y)}{\lg n} p^n(y-z) (1+ o(1)). $$
Since  $p^n(y-z) = p^n(y-x)(1+o(1))$ if  both sides are positive, 
we may write (\ref{13}) as 
$$p_A^n(x,y) = P_x[\tau_{ U(R)} < N\wedge \sigma_A]\frac{2\k \;\! u_{-A}(-y)}{\lg n} p^n(y-x)(1+ o(1)) + \e(n,x,y).  $$
 The evaluation of the term  $\e(n,x,y) $ given in (\ref{14})   and that of
  $P_x[\tau_{ U(R)}  < N\wedge\sigma_A]$   in (\ref{24}) are both available. Thus   we conclude the required relation of  Theorem \ref{thm1}. \qed

\v2\v2
{\bf Acknowledgments.}\, 
 I wish to thank  an anonymous  referee for his/her   helpful comments.


\v2


\begin{thebibliography}{99}
\baselineskip=9pt










\bibitem {F} W. Feller, Introduction to probability theory and its applications, vol. II, 2nd ed., John Wiley \& Sons, New York 1971.


\bibitem {K} H. Kesten,  Ratio limit theorems II, Journal d'Analyse Math. {\bf 11}, (1963), 323-379.

\bibitem {KSK} J. Kemeny, J. Snell and A. Knapp, denumerable Markov chains, Van Nostrand Co., Princeton, 1966.



\bibitem {Ht}  G. Hunt, Some theorems concerning Brownian motion, TAMS. {\bf s1} (1956),  294-319.



\bibitem {LL} G. F. Lawler and V. Limic, Random walk: a modern introduction, Cambridge Univ. Press,  2010. 






\bibitem {Sh} F. Spitzer, Hitting probabilities, J. Math. and Mech., {\bf 11} (162), 593-614.

\bibitem {S} F. Spitzer,  Principles of Random Walk, Van Nostrand, Princeton, 1964. 



\bibitem {St} D. W. Stroock, Probability, 2nd ed., Cambridge Univ. Press, 2011.







\bibitem{Ugf}  K. Uchiyama, Green's functions for random walks on $\Z^d$, Proc. London Math. Soc. {\bf 77} (1998), 215-240.

\bibitem {Um_add}  K. Uchiyama,, Asymptotic estimates of the Green functions and transition
 probabilities for Markov additive processes,  Electric Journal
of Probability {\bf 12} 
(2007), 138-180


\bibitem{Ufh}  K. Uchiyama, The first  hitting time of a single point for random walks. Elect. J. Probab. {\bf 16}, 
(2011), 1160-2000 

\bibitem{Ubdsty} K. Uchiyama, The   transition density  of Brownian  motion  killed on  a bounded set, preprint, available at: http://arxiv.org/abs/1603.03902

 

\bibitem{U1dm}  K. Uchiyama, One dimensional  random walk  killed on a finite set, preprint,
available at: http://arxiv.org/abs/1603.02117
















%




\end{thebibliography}
 \end{document}